\documentclass[12pt]{article}
\usepackage{e-jc}

\usepackage{amsmath,url,color,tikz}
\usepackage{graphicx}
\usepackage{lmodern}
\usepackage{epstopdf}
\usepackage{array}
\usepackage{mathtools}
\usepackage{enumerate}
\usepackage{lscape}
\usepackage{verbatim}
\usepackage{multicol}

\dateline{TBD}{TBD}{TBD}


\Copyright{Jesse Kim and James Propp. Released under the CC BY-ND license (International 4.0).}

\title{A pentagonal number theorem for tribone tilings}

\author{Jesse Kim
\\
\small Department of Mathematics \\[-0.8ex]
\small University of California, San Diego\\[-0.8ex]
\small California, U.S.A.\\
\small\tt jvkim@ucsd.edu \\
\and
James Propp\thanks{Supported by a Simons Collaboration Grant.}\\
\small Department of Mathematical Sciences\\[-0.8ex]
\small University of Massachusetts, Lowell\\[-0.8ex]
\small Massachusetts, U.S.A.\\
\small\tt jamespropp@gmail.com}

\begin{document}

\maketitle


\begin{abstract}
\noindent
Conway and Lagarias showed that certain roughly triangular regions 
in the hexagonal grid cannot be tiled by shapes Thurston later dubbed tribones.
Here we study a two-parameter family of roughly hexagonal regions 
in the hexagonal grid and show that a tiling by tribones exists
if and only if the two parameters associated with the region are 
the paired pentagonal numbers $k(3k \pm 1)/2$.
\end{abstract}

\newcommand{\RR}{\mathbb{R}}
\newcommand{\ZZ}{\mathbb{Z}}
\newcommand{\CC}{\mathbb{C}}
\newcommand{\ba}{{\rm \bf{a}}}
\newcommand{\bb}{{\rm \bf{b}}}
\newcommand{\bc}{{\rm \bf{c}}}
\newcommand{\bi}{{\rm \bf{i}}}
\newcommand{\bj}{{\rm \bf{j}}}
\newcommand{\bv}{{\rm \bf{v}}}
\newcommand{\bw}{{\rm \bf{w}}}
\newcommand{\bzero}{\rm \bf{0}}
\hyphenation{hexa-gonal}
\hyphenation{Thur-ston}
\hyphenation{re-pre-sen-ta-tive}







\section{Introduction} \label{sec-intro}

\begin{figure}[h]
\begin{center}
\includegraphics[width=2.0in]{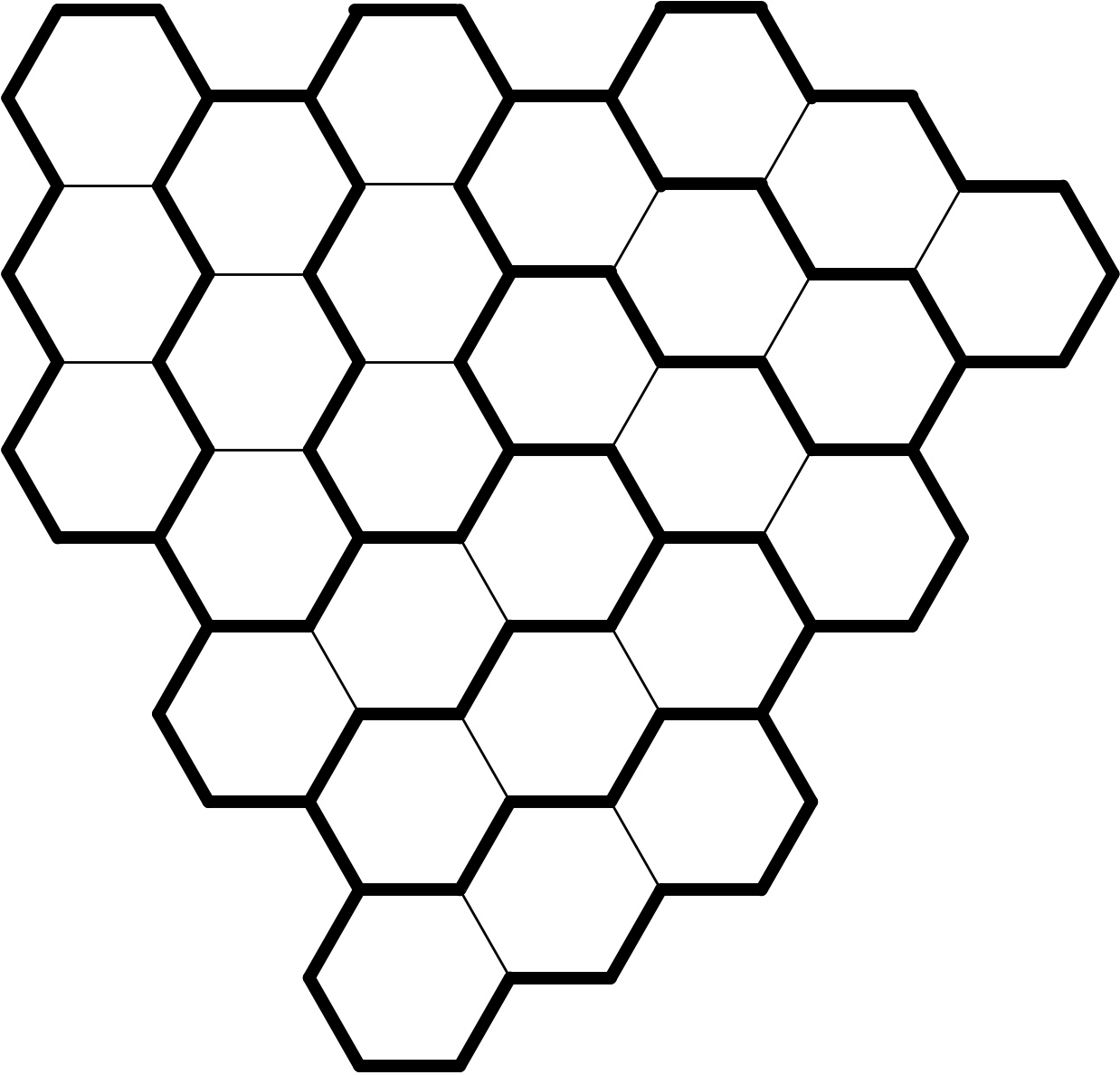}
\end{center}
\caption{The (5,7)-benzel tiled by tribones
in one of the two possible ways.}
\label{fig:what-is}
\end{figure}

\begin{figure}[h]
\begin{center}
\includegraphics[width=5.0in]{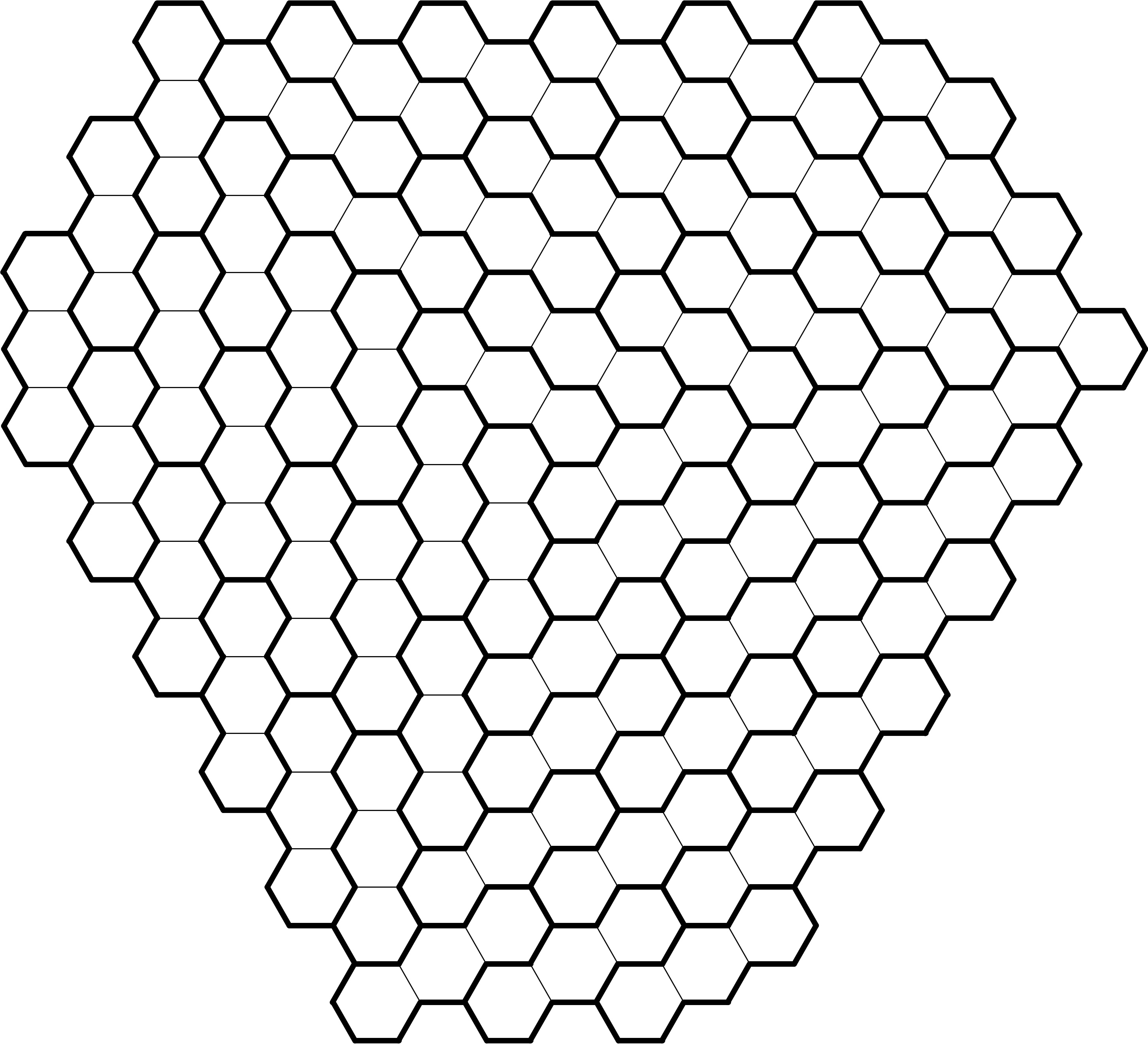}
\end{center}
\caption{The (12,15)-benzel tiled by tribones
in one of the 42,705 possible ways.}
\label{fig:twelve-fifteen}
\end{figure}

\begin{figure}
\begin{center}
\includegraphics[width=2.0in]{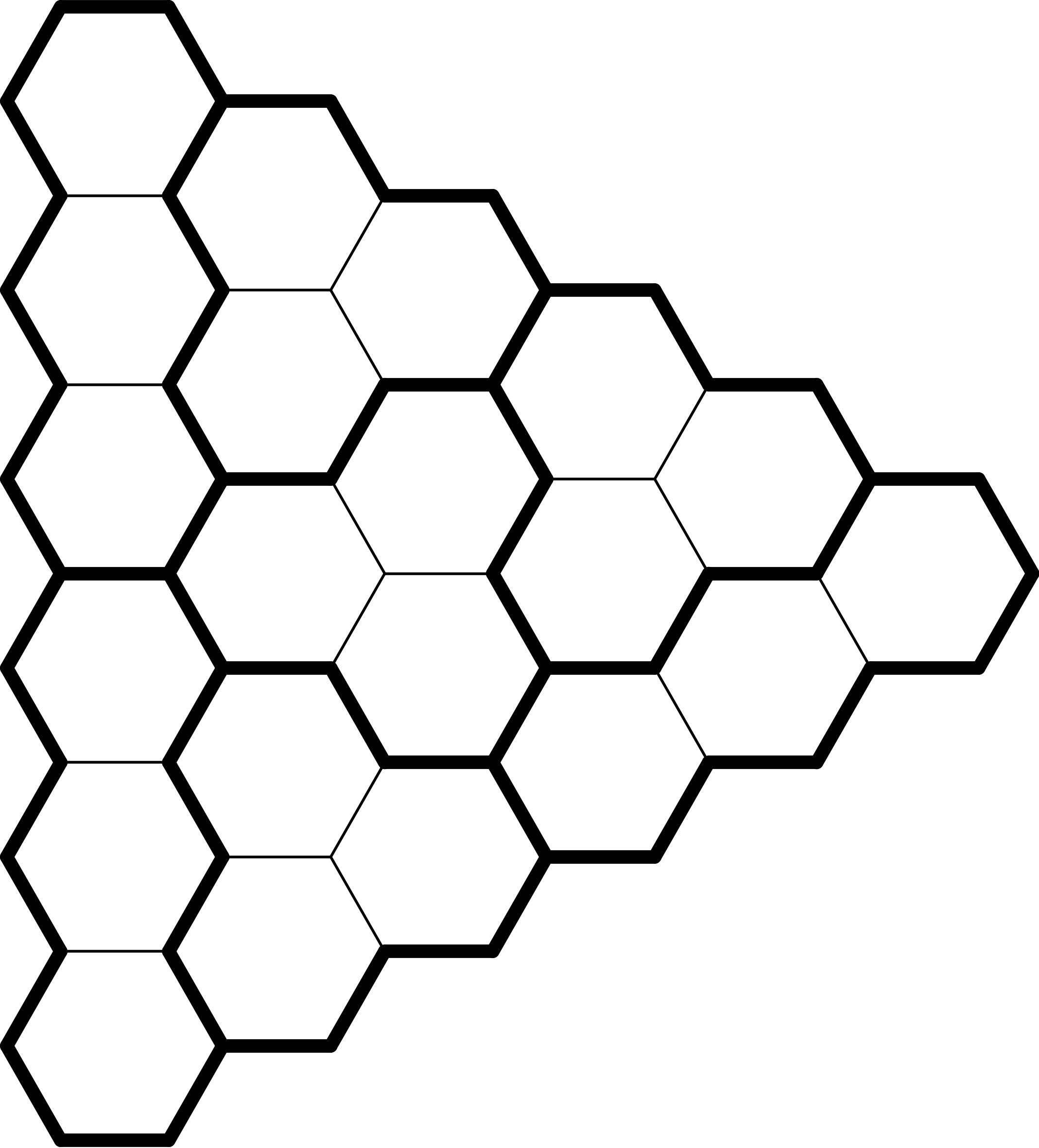}
\end{center}
\caption{The honeycomb triangle $T_6$ tiled by stones and bones.}
\label{fig:tee-six}
\end{figure}

This article treats tilings such as the one shown
in Figures~\ref{fig:what-is} and~\ref{fig:twelve-fifteen}
in which the tiles are {\em tribones}
consisting of three successive unit hexagons with collinear centers.
(Here we use ``tribone'' in the sense introduced by Thurston~\cite{Thur};
we will often call such tiles simply {\em bones}.)
Tribone tilings were studied by Conway and Lagarias
who showed that for all integers $N \geq 1$,
the roughly triangular region $T_N$
exemplified for $N=6$ in Figure~\ref{fig:tee-six}
cannot be tiled by tribones (which they called $L_3$'s).
In Figure~\ref{fig:tee-six} the region $T_6$ is tiled by a mixture
of tribones and additional tiles Conway and Lagarias called $T_2$'s
(which we will call {\em stones}),
where three of the stones point to the right and one points to the left.
Conway and Lagarias showed that 
the number of right-pointing stones minus the number of left-pointing stones
is the same for all tilings of a simply-connected region.
Since this quantity (a rescaled version of the {\em Conway-Lagarias invariant})
is nonzero for the particular tiling of $T_6$ shown in Figure~\ref{fig:tee-six},
it must be nonzero for all tilings of $T_6$,
and in particular, there can be no pure tribone tilings.  Similar arguments
show that $T_N$ cannot be tiled by tribones for any positive value of $N$.

(In the preceding paragraph and throughout this paper 
we conflate the subtly different perspectives
of Conway and Lagarias~\cite{CoLa} and Thurston~\cite{Thur}
while introducing some notions of our own;
readers who want a fine-grained understanding
of the lineage of the ideas in this paper
should read those two articles, as well as the article~\cite{Pro1}.)

Here we study a two-parameter family of regions
introduced in~\cite{Pro2} called {\em benzels};
for instance, Figure~\ref{fig:what-is} shows the (5,7)-benzel
and Figure~\ref{fig:twelve-fifteen} shows the (12,15)-benzel.
A variant of Gauss' shoelace formula
allowed the first author one to compute the signed area (aka algebraic area)
enclosed by a closed polygonal path
and, by ``twisting'' the formula in the manner prescribed by
the work of Conway, Lagarias, and Thurston,
to compute the values of the Conway-Lagarias invariant for all benzels.
This enabled him to prove in the spring of 2022
that the $(a,b)$-benzel cannot be tiled by tribones
unless $\{a,b\} = \{k(3k+1)/2,k(3k-1)/2\}$ for some positive integer $k$,
that is, unless $a$ and $b$ are paired pentagonal numbers.
It remained to show that when $a$ and $b$ are of this form,
at least one tribone tiling exists; several students at the
Open Problems in Algebraic Combinatorics 2022 meeting rose to the challenge.
The first author's construction was the clearest,
and he also provided an argument explaining why the construction always works.
These two results together yielded our main result (Theorem 3 below):

\begin{center}
\fbox{
\parbox{250pt}{
The $(a,b)$-benzel can be tiled by tribones if and \\
only if $a$ and $b$ are paired pentagonal numbers.
}}
\end{center}

\noindent
In referring to this result as a ``pentagonal number theorem''
we are hearkening back to Euler's famed theorem
about the series expansion of $\prod_{n=1}^{\infty} (1-q^n)$,
but we know of no connection between the two theorems.

Section~\ref{sec-benzels} defines benzels,
which are parametrized by pairs $a,b$ 
where $a$ and $b$ are positive integers
satisfying $2 \leq a \leq 2b$ and $2 \leq b \leq 2a$.
In section~\ref{sec-area} we express 
the boundary of the $(a,b)$-benzel
via a concatenation of unit vectors
and compute the signed area enclosed by the boundary path
(Theorem 1).
Section~\ref{sec-shadow} presents a reinterpretation 
of the Conway-Lagarias invariant
as the signed area enclosed by the ``shadow path''
associated with the boundary of the region being tiled.
Section~\ref{sec-formula} combines the results 
of the three preceding sections to obtain 
a formula for the value of the Conway-Lagarias invariant
(Theorem 2).

Section~\ref{sec-vanish} analyzes the polynomials
computed in section~\ref{sec-formula} 
to show that the Conway-Lagarias invariant vanishes
only when $a=k(3k-1)/2$ and $b=k(3k+1)/2$ (or vice versa) for some $k$.
Section~\ref{sec-construct} shows, conversely,
that at least one tribone tiling exists
when $a$ and $b$ are of this form.
This proves that the algebraic condition on $a$ and $b$
is both necessary and sufficient for the existence
of tilings of the $(a,b)$-benzel by bones (Theorem 3).

The companion article~\cite{Pro2} discusses 
other classes of stones-and-bones tilings of benzels,
offering some (mostly conjectural) exact formulas for special cases 
in which certain prototiles are permitted/prohibited
and $a$ and $b$ take on certain kinds of values.
The two articles adopt slightly different conventions.
In that other article, the focus is on the hexagonal cells
and it makes sense to locate the three middle cells of a benzel
so that their centers are at 1, $\omega$, and $\omega^2$ in the complex plane.
In this article the vertices and edges along the boundaries of the cells
play a leading role so it makes more sense to have edges joining 
the vertex 0 to the vertices 1, $\omega$, and $\omega^2$ 
in the benzel graph dual to the benzel 2-complex.
The reader should have little difficulty
translating between the two pictures
using appropriate reflections and sign-changes where needed.

\section{Stones, bones, and benzels} \label{sec-benzels}

In this section we define the sorts of regions we seek to tile
and the sorts of tiles we are allowed to use.

\begin{figure}
\begin{center}
\includegraphics[width=4.0in]{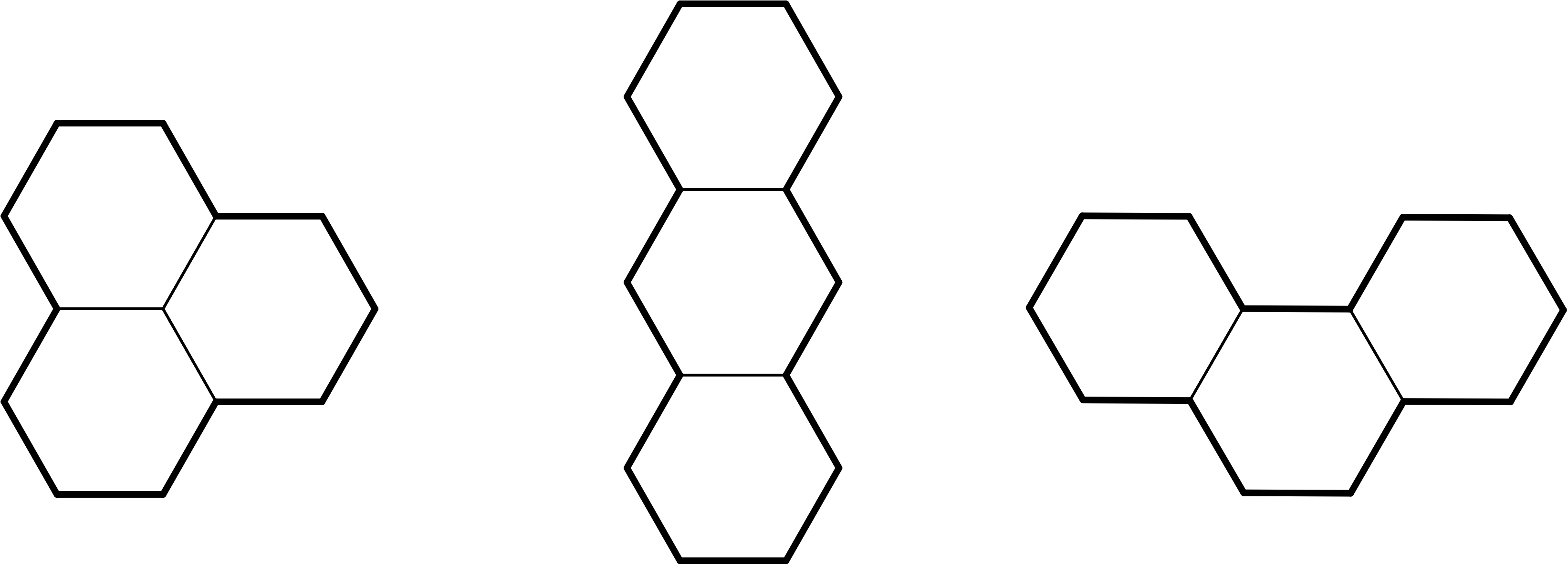}
\end{center}
\caption{The three kinds of trihex: stone, bone, and phone.} 
\label{fig:trihex}
\end{figure}

Conway and Lagarias~\cite{CoLa} refer to the first two tiles 
shown in Figure~\ref{fig:trihex} as $T_2$ and $L_3$.
Thurston~\cite{Thur} keeps the name $T_2$ but calls the $L_3$ a tribone.
Ardila and Stanley~\cite{ArSt} call the $T_2$ a tribone and the $L_3$ a trihex.
Meanwhile, in the recreational mathematics literature,
a trihex denotes any of the three tiles shown in Figure~\ref{fig:trihex}.
Curiously, although the literature on polyforms has given names 
to the tetrahexes, it has not given names to the trihexes. 
The second author has dubbed them the {\em stone}, the {\em bone}, 
and the {\em phone}, respectively~\cite{Pro1}.
Phones will not play a role in this article.
``Bone'' is sometimes used as a synonym for domino
but we think this ambiguity is unlikely to cause confusion
since the contexts are very different.
We define the {\em axis} of a bone as
the line through the centers of its three constituent hexagons.

\begin{figure}
\begin{center}
\includegraphics[width=2.5in]{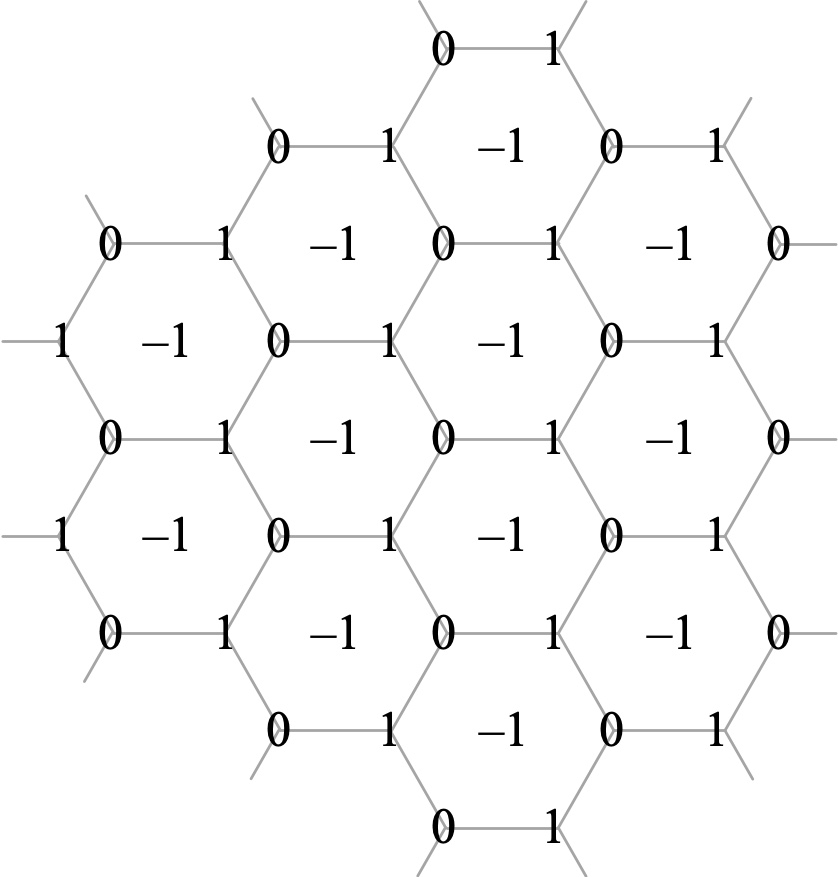}
\end{center}
\caption{The lattice $L=\ZZ + \ZZ\omega$.}
\label{fig:sublattices}
\end{figure}

We embed our tilings in the complex plane.  Let $\omega = e^{2 \pi i/3}$;
all points of interest will belong to the lattice $L = \ZZ+\ZZ\omega$.
More specifically, consider the index-three sublattice 
$L_0 = \ZZ(1-\omega)+\ZZ(1-\omega^2)$
which together with its translates 
$L_1 = L_0+1$ and $L_{-1} = L_0-1$ partition $L$;
we form a graph whose vertices are the elements of $L_0 \cup L_1$ 
with edges connecting two vertices whenever the distance between them is 1.
A point $\alpha$ in $L_{-1}$ is the center of a hexagon with corner points 
$\alpha \pm 1$, $\alpha \pm \omega$, and $\alpha \pm \omega^2$;
a point $\alpha$ in $L_0$ borders three hexagons, 
one of which is to its left;
and a point $\alpha$ in $L_1$ borders three hexagons, 
one of which is to its right.
See Figure~\ref{fig:sublattices},
in which points are marked 0, $1$, or $-1$
according to whether they belong to $L_0$, $L_1$, or $L_{-1}$.

\begin{figure}
\begin{center}
\includegraphics[width=4.0in]{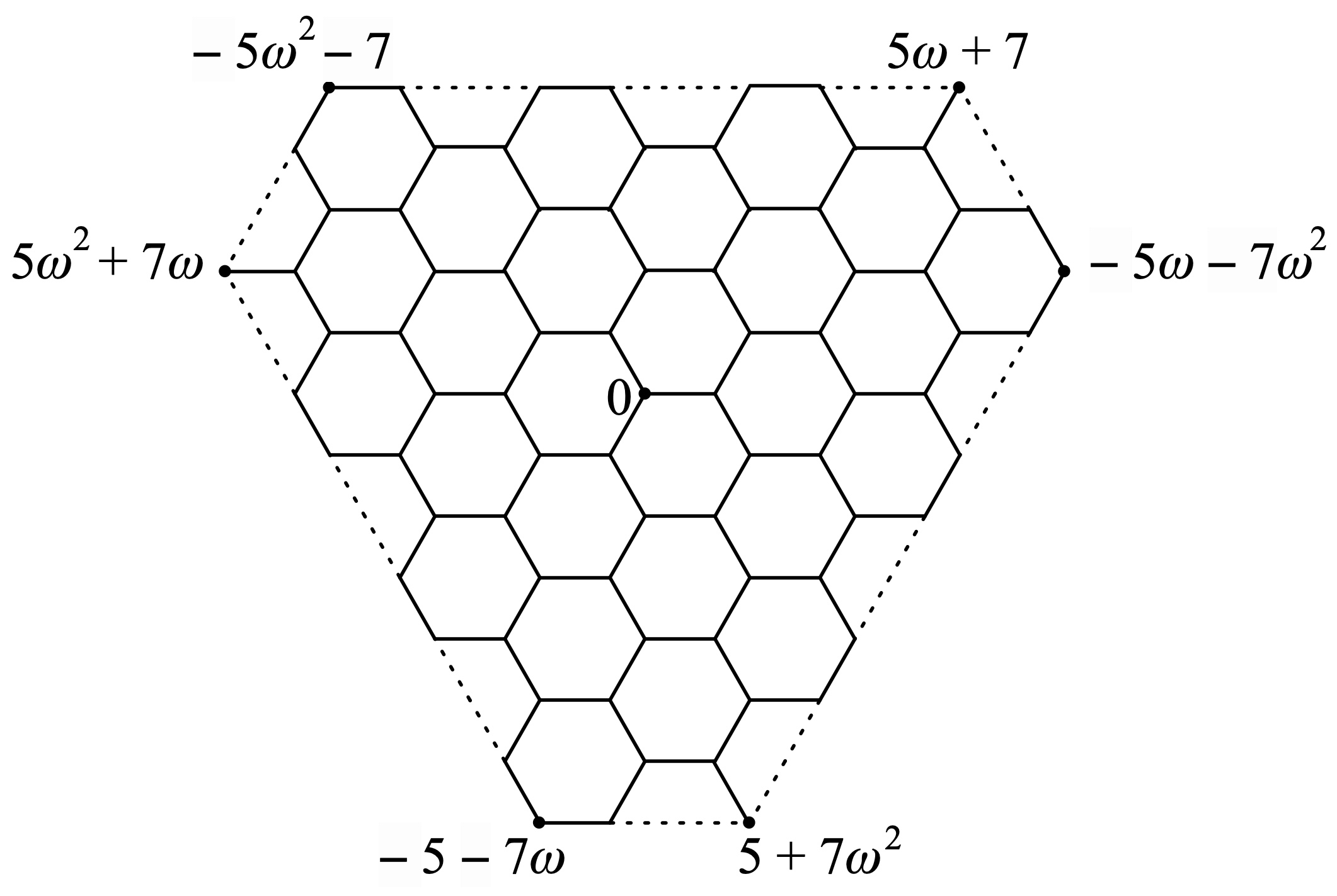}
\end{center}
\caption{The (5,7)-benzel and its bounding hexagon.}
\label{fig:corners}
\end{figure}

For positive integers $a,b$ satisfying
$2 \leq a \leq 2b$ and $2 \leq b \leq 2a$
we define the $(a,b)$-benzel as the union of all the unit hexagons
that lie completely within the large hexagon with vertices
$a\omega+b$, $-a\omega^2-b$,
$a\omega^2+b\omega$, $-a-b\omega$,
$a+b\omega^2$, and $-a\omega-b\omega^2$,
as shown in Figure~\ref{fig:corners} for $a=5$, $b=7$.
We refer to this hexagon as the {\em bounding hexagon} of the benzel.
The bounding hexagon has three sides of length $2a-b$
alternating with three sides of length $2b-a$.
Note that the figure includes three stray edges, or ``spurs'',
at $5\omega+7$, $5\omega^2+7\omega$, and $5+7\omega^2$;
although such edges are not part of the proper boundary of the benzel,
it is sometimes convenient to include them.

We say that a benzel is of class 0, 1, or $-1$
according to the residue of $a+b$ mod 3.
If a benzel is of class $c$,
then three of the vertices of the bounding hexagon
($a\omega+b$, $a\omega^2+b\omega$, and $a+b\omega^2$) belong to $L_c$
and the other three 
($-a\omega^2-b$, $-a-b\omega$, and $-a\omega-b\omega^2$) belong to $L_{-c}$.
Note that $2a-b \equiv 2b-a \equiv -c$ mod 3.

\begin{figure}
\begin{center}
\includegraphics[width=5.0in]{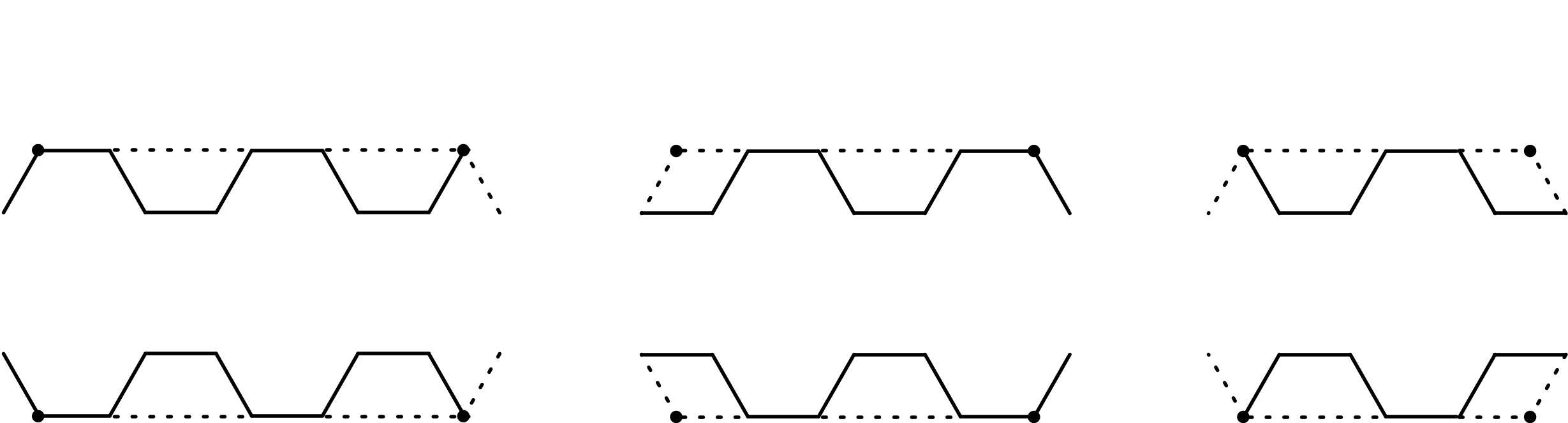}
\end{center}
\caption{The upper and lower stretches of the boundary of a class 0 benzel (left),
a class 1 benzel (middle), and a class $-1$ benzel (right).}
\label{fig:three-cases}
\end{figure}

The $(a,b)$-benzel has 120-degree rotational symmetry.
The boundary of a benzel consists of six stretches
that follow the successive sides of the bounding hexagon.
The six panels of Figure~\ref{fig:three-cases}
show the upper and lower stretches of the boundary of the $(a,b)$-benzel 
in the three cases $c=0$, $c=1$, and $c=-1$. 
The upper and lower boundaries of the benzel
contain $(2b-a+c)/3$ and $(2a-b+c)/3$ horizontal edges respectively.
The other four stretches are just rotated copies of those two;
the six stretches alternate between
ones that look like the upper stretch
and ones that look like the lower stretch.
In these examples we have taken
$a=b=6$, $a=b=5$, and $a=b=4$, respectively,
resulting in upper and lower stretches of the same length;
when $a \neq b$ the upper and lower stretches 
will not have the same length
but will still exhibit the behavior seen in the figure
with regard to the appearance
of the left and right ends of the upper and lower stretches.

\section{A Gauss formula for the hexagonal grid} \label{sec-area}

Our goal in this section is to derive
a formula for the signed area enclosed
by the counterclockwise-oriented boundary of a benzel
(see Theorem 1 at the end of this section).

For regions described by boundary-paths made of short segments,
it's helpful to express the signed area as a double-sum of
scalar cross-products of those segments
(where the scalar cross-product $\bv_1 \times \bv_2$ of 
$\bv_1 = (x_1,y_1)$ and $\bv_2 = (x_2,y_2)$ in $\RR^2$
is defined as $x_1 y_2 - x_2 y_1$).
To be specific, suppose we have vectors $\bv_i$ ($1 \leq i \leq n$) satisfying
$\sum_{i=1}^n \bv_i = \bzero$. Then, with $\bw_i = \bv_1 + \dots + \bv_i$,
the signed area enclosed by the $n$-gon with
vertices $\bw_1$, $\bw_2$, \dots, $\bw_n = \bzero$
is $\frac12 \sum_{i} \bw_i \times \bv_{i+1}$, which can be rewritten as
$\frac12 \sum_{1 \leq i < j \leq n} \bv_i \times \bv_j$.
One may call this (with slight historical distortion)
a {\em Gauss area formula} for the signed area enclosed by the closed curve.
Note that if the closed curve is simple (i.e., does not intersect itself)
and encircles a region in the counterclockwise direction,
the signed area associated with the curve
is the ordinary unsigned area of the region enclosed.

As an example, consider the square grid with unit vectors
$\bi = (1,0)$ and $\bj = (0,1)$, and consider the $k$-by-$k$ square
with vertices $(0,0)$, $(k,0)$, $(k,k)$, $(0,k)$, and $(0,0)$.
The counterclockwise boundary of the square consists of
$k$ occurrences of $\bi$, $k$ occurrences of $\bj$,
$k$ occurrences of $-\bi$, and $k$ occurrences of $-\bj$, in that order.
If we ignore the terms of the double sum
$\sum_{1 \leq i < j \leq n} \bv_i \times \bv_j$
that vanish because they are of the form $\bv \times \bv$
or $\bv \times -\bv$,
we obtain $k^2$ occurrences of $\bi \times \bj$,
$k^2$ occurrences of $\bi \times -\bj$,
$k^2$ occurrences of $\bj \times -\bi$,
and $k^2$ occurrences of $-\bi \times -\bj$.
This amounts to $k^2$ times $1 - 1 + 1 + 1 = 2$,
so we find that the signed area is $\frac12 2k^2 = k^2$, as it should be.
For brevity, we will write the boundary as
$(\bi)^k (\bj)^k (-\bi)^k (-\bj)^k$,
where superscripts denote string repetition.
A more complicated example is the Aztec diamond~\cite{EKLP} of order $k$;
the reader may wish to check that applying the Gauss formula to
the string $(\bi,\bj)^k (\bj,-\bi)^k (-\bi,-\bj)^k (-\bj,\bi)^k$
gives the predicted signed area $2k^2+2k$.
For a general path whose segments are $\pm \bi$ and $\pm \bj$ ,
the signed area $\frac12 \sum_{1 \leq i < j \leq n} \bv_i \times \bv_j$
is equal to the sum, over all square cells $C$,
of the winding number of the path around $C$;
although there are infinitely many cells $C$ in the infinite square grid,
the winding number of any particular finite path around $C$
vanishes for all but finitely many $C$, so the sum is well-defined.

\begin{figure}
\begin{center}
\includegraphics[width=2.0in]{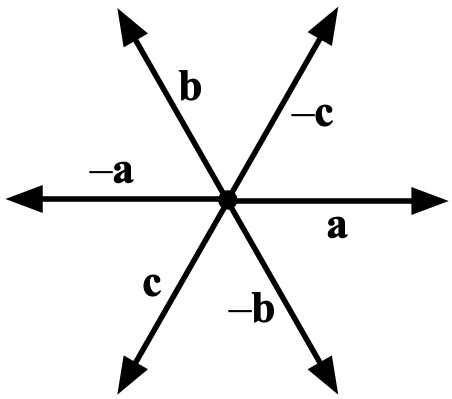}
\end{center}
\caption{The six unit vectors.}
\label{fig:six-vectors}
\end{figure}

\begin{figure}
\begin{center}
\includegraphics[width=2.0in]{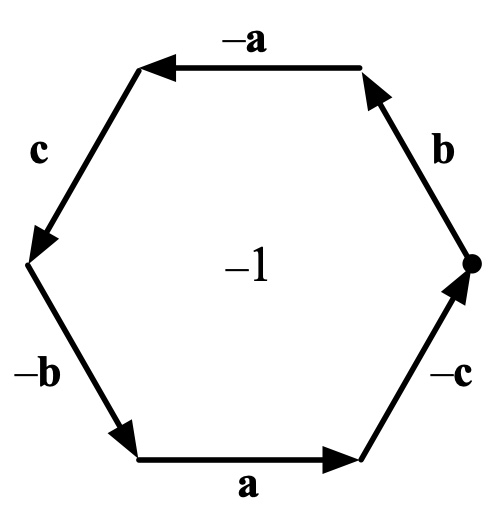}
\end{center}
\caption{The counterclockwise path enclosing a hexagon.}
\label{fig:one-hexagon}
\end{figure}

For the hexagonal grid we will need 
three unit vectors $\ba$, $\bb$, and $\bc$,
pointing from $0$ to $1$, $\omega$ and $\omega^2$,
respectively, under our identification of $\CC$
with the plane containing the tiling;
see Figure~\ref{fig:six-vectors}.
To remove factors of $\sqrt{3}$,
we redefine $\times$ so that
$\ba \times \bb = \bb \times \bc = \bc \times \ba = +1$,
$\bb \times \ba = \bc \times \bb = \ba \times \bc = -1$,
and $\ba \times \ba = \bb \times \bb = \bc \times \bc = 0$
(with other cases being determined by the assumption
that $\times$ is bilinear); this bilinear form is antisymmetric.
If each hexagonal cell is assigned area 1,
then the signed area enclosed by a path
given by the vectors $\bv_1,\dots,\bv_n$
(with $\bv_i \in \{\pm \ba, \pm \bb, \pm \bc\}$ for $1 \leq i \leq n$)
is $\frac16 \sum_{1 \leq i < j \leq n} \bv_i \times \bv_j$.
This is equal to the sum, over all hexagonal cells $C$,
of the winding number of the path around $C$.
(Here we are assuming that the path
is the counterclockwise boundary of a union of hexagons.
For instance, we are not considering the closed path $\ba,\bb,\bc$
which encircles a triangle.)
As an example, consider the unit hexagon centered at $-1$
as in Figure~\ref{fig:one-hexagon}.
The counterclockwise path around the hexagon
that starts and ends at the origin
(the rightmost corner of the hexagon)
is associated with the sequence of vectors $\bb,-\ba,\bc,-\bb,\ba,-\bc$.
The nonvanishing terms of the sum
$\sum_{1 \leq i < j \leq 6} \bv_i \times \bv_j$ are 
$\bb \times -\ba$, $\bb \times \bc$, $\bb \times \ba$,
$\bb \times -\bc$, $-\ba \times \bc$, $-\ba \times -\bb$,
$-\ba \times -\bc$, $\bc \times -\bb$, $\bc \times \ba$,
$-\bb \times \ba$, $-\bb \times -\bc$, and $\ba \times -\bc$.
Using the fact that $\times$ is antisymmetric,
we can simplify the sum, obtaining
$2( \ba \times \bb + \bb \times \bc + \bc \times \ba) = 6$,
so that $\frac16 \sum_{1 \leq i < j \leq n} \bv_i \times \bv_j = 1$,
which is indeed the sum of the winding numbers
of the paths around all the hexagonal cells,
which is 1 for the cell centered at $-1$ and 0 for all other cells.

To save space on the page,
we will sometimes write $-\ba$, $-\bb$, and $-\bc$
as $\ba'$, $\bb'$, and $\bc'$, respectively.

Take $a,b$ with $a+b \equiv 0$ mod 3.
The upper stretch of the $(a,b)$-benzel,
starting at the upper right corner of the bounding hexagon
and ending at the upper left corner of the bounding hexagon
(see the first panel of Figure~\ref{fig:three-cases})
is $(\bc,\ba',\bb,\ba')^t$ with $t = (2b-a)/3$.
The lower stretch of the $(a,b)$-benzel,
starting at the lower left corner of the bounding hexagon
and ending at the lower right corner of the bounding hexagon
is $(\ba,\bc',\ba,\bb')^s$ with $s = (2a-b)/3$.
Thus the counterclockwise boundary of the benzel,
starting and ending at the rightmost corner of the bounding hexagon,
is the concatenation
$$(\bb,\ba',\bb,\bc')^s (\bc,\ba',\bb,\ba')^t
  (\bc,\bb',\bc,\ba')^s (\ba,\bb',\bc,\bb')^t
  (\ba,\bc',\ba,\bb')^s (\bb,\bc',\ba,\bc')^t$$
involving $12s+12t$ unit vectors.
One could in principle work out the value of the double sum
$\frac16 \sum_{1 \leq i < j \leq 12s+12t} \bv_i \times \bv_j$,
but it is also possible to determine the general form of the answer
and then fit the undetermined coefficients.
Specifically, one can lump together 
the ${12s+12t \choose 2}$ products $\bv_i \times \bv_j$ into 
${6 \choose 2} \cdot 4 \cdot 4 + {6 \choose 1} \cdot 4 \cdot 4 = 256$ 
types (according to which of the 24 symbols in the above expression
give rise to $\bv_i$ and $\bv_j$ respectively)
such that the number of terms $\bv_i \times \bv_j$ of each type
is a quadratic function of $s$ and $t$.
This implies that the area of the benzel
is itself a quadratic function of $s$ and $t$,
and hence a quadratic function of $a$ and $b$.
This quadratic function has six undetermined coefficients
that can be determined 
by finding the area of half a dozen specific benzels
and solving the resulting system of linear equations.
One finds in this way that the area of a class 0 benzel 
is equal to $(-a^2+4ab-b^2-a-b)/2$
(a positive integer divisible by 3).

One can analyze benzels of class 1 and class $-1$ in a similar way.
For both of these cases, one can augment the boundary path
by spurs that join the path to one of the corners of the bounding hexagon;
these spurs are traversed in one direction
and then immediately traversed in the opposite direction,
so they enclose no cells and make no contribution to the area.

\bigskip

{\bf Theorem 1}: Let $R$ be a benzel of class $c$.  Then the area of $R$ is 
$$\left\{ \begin{array}{ll}
(-a^2+4ab-b^2-a-b)/2 & \mbox{if $c$ is 0 or $-1$}, \\
(-a^2+4ab-b^2-a-b+2)/2 & \mbox{if $c$ is 1.}
\end{array} \right.$$

\section{Shadow words and the Conway-Lagarias invariant} \label{sec-shadow}

In this section we will lay the groundwork for computing 
the Conway-Lagarias invariant of benzels in the next section.

Conway and Lagarias show that in any simply-connected region $R$
that can be tiled by stones and bones,
the number of right-pointing stones
minus the number of left-pointing stones
depends only on the region being tiled,
not the particular tiling.  We call this the 
{\em rescaled Conway-Thurston invariant} and denote it by $i(R)$.
Thurston shows how to compute this invariant
using the Cayley graphs of the two groups
$$A = \langle a,b,c \: | \: a^2 = b^2 = c^2 = (abc)^2 = 1 \rangle$$
and $$T_0 = \langle a,b,c \: | \: a^2 = b^2 = c^2 
= (ab)^3 = (bc)^3 = (ca)^3 = 1 \rangle$$
(where the formal involutions $a$, $b$, and $c$ 
are not to be confused with the integers $a,b,c$ 
introduced in section~\ref{sec-benzels}
but are closely related to the lax vectors
$\pm \ba$, $\pm \bb$, and $\pm \bc$ with $\ba$, $\bb$, and $\bc$
as in section~\ref{sec-area}; see~\cite{Con} for more on lax vectors).
Specifically, one draws $R$ in $\Gamma(A)$, the Cayley graph of $A$,
and picks a basepoint $p$ on the boundary of $R$;
one expresses the boundary path that goes around $R$
once counterclockwise from $p$ to itself
as a product of $a$'s, $b$'s, and $c$'s;
one reinterprets the letters of this word as the generators of the group $T_0$;
one draws a path in $\Gamma(T_0)$, the Cayley graph of $T_0$, 
associated with that word
(the new path must in fact be a closed path, for non-obvious reasons);
and lastly one computes the signed area enclosed by the new path;
this signed area, which we call the {\em (unrescaled) Conway-Lagarias invariant}
and denote by $I(R)$, is equal to $3 i(R)$.

\begin{figure}
\begin{center}
\includegraphics[width=4.0in]{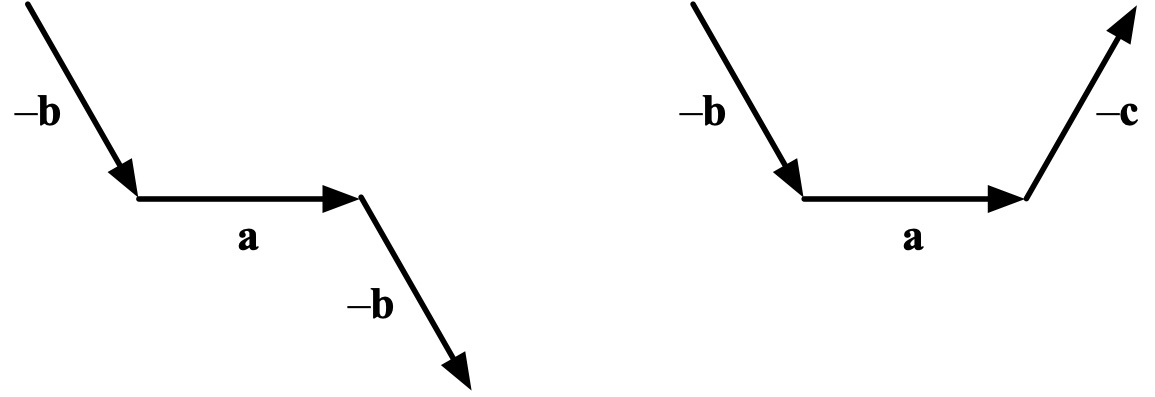}
\end{center}
\caption{Weaving versus winding.}
\label{fig:weave-wind}
\end{figure}

Thurston remarks ``It is curious that even though the groups $A$ and $T_0$
and the labeled graphs $\Gamma(A)$ and $\Gamma(T_0)$ are different,
when the labels are stripped they become isomorphic.''
Specifically, both are isomorphic to the hexagon graph,
and there is a direct pictorial way to convert the old path into the new
without explicit recourse to the groups $A$ and $T_0$.
Given a basepoint $p \in L_0$ on the boundary of $R$,
let $e_1,\dots,e_n$ be the sequence of edges
traversing the boundary of $R$ counterclockwise
from $p$ back to itself, and let $e_0=e_n$ and $e_{n+1}=e_1$.
Call this path $\pi$.
For $1 \leq i \leq n$, say that $\pi$ {\em weaves} at step $i$
if $e_{i-1}$ and $e_{i+1}$ are parallel
and say that $\pi$ {\em winds} otherwise,
i.e., if $e_{i-1}$, $e_i$, and $e_{i+1}$
are consecutive edges of some hexagon;
see Figure~\ref{fig:weave-wind},
in which the edges are given orientations
and marked with the associated vectors.
Suppose we have another path
$e'_0,\dots,e'_n$ (call it $\pi'$ for short)
in the hexagon graph that starts from a vertex $p' \in L_0$,
and suppose that $\pi'$ {\em shadows} $\pi$
for all $1 \leq i \leq n$ in the sense that
$\pi'$ weaves at step $i$ if $\pi$ winds at step $i$
and $\pi'$ winds at step $i$ if $\pi$ weaves at step $i$.
(This is the ``twist'' referred to in section~\ref{sec-intro}.)
Then the unrescaled Conway-Lagarias invariant of $R$
equals the signed area enclosed by $\pi'$.

\begin{figure}
\begin{center}
\includegraphics[width=4.0in]{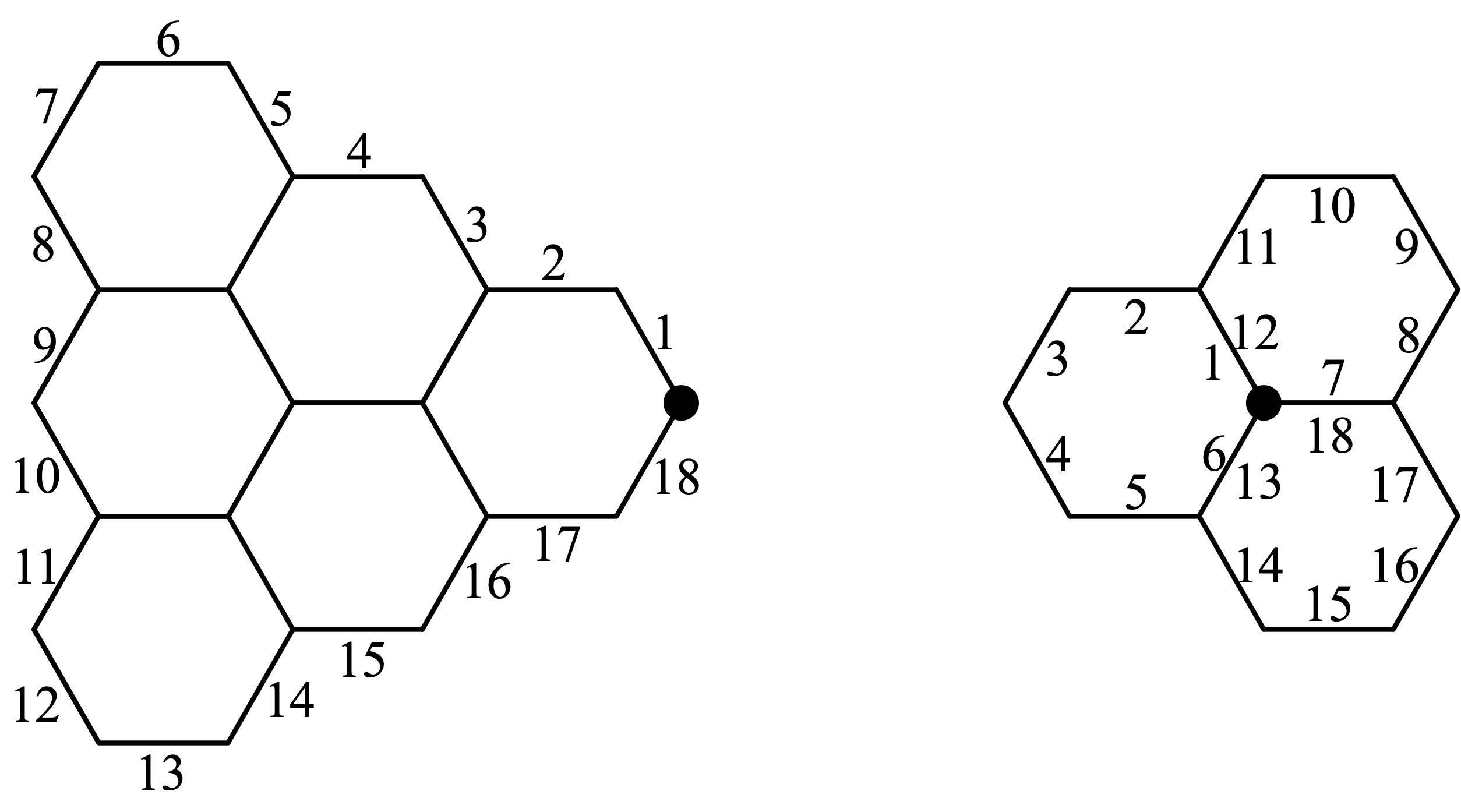}
\end{center}
\caption{A path and its shadow.}
\label{fig:tee-three}
\end{figure}

For example, suppose $R$ is the (3,3)-benzel (aka $T_3$),
as shown in the left panel of Figure~\ref{fig:tee-three},
let the basepoint $p$ be its rightmost point,
and let $\pi$ be the path that encircles $R$ counterclockwise
from $p$ to itself, with edges numbered 1 through 18.
Let $\pi'$ be the path shown in the right panel of Figure~\ref{fig:tee-three},
starting and ending at the marked point.
It is easy to check that $\pi'$ shadows $\pi$ and that $\pi'$ encircles
three hexagons in the counterclockwise direction, enclosing signed area 3.
This is also the total area of all the right-pointing stones
minus the total area of all the left-pointing stones
in each of the (three) tilings of the (3,3)-benzel by stones and bones.
Note that in terms of unit vectors, the path $\pi$ is given by
$$\bb,-\ba,\bb,-\ba,\bb,-\ba,\bc,-\bb,\bc,-\bb,\bc,-\bb,\ba,-\bc,\ba,-\bc,\ba,-\bc\:$$
while the path $\pi'$ is given by
$$\bb,-\ba,\bc,-\bb,\ba,-\bc,\ba,-\bc,\bb,-\ba,\bc,-\bb,\bc,-\bb,\ba,-\bc,\bb,-\ba;$$
such representations of paths as strings of unit vectors will be helpful shortly.

It should be stressed that $p$ can be any point in $\partial(R) \cap L_0$
and $p'$ can be any point in $L_0$.
Moreover, for fixed $p$ and $p'$,
there are six paths $\pi'$ that shadow a given $\pi$:
if $e'_1$ is any of the three edges emanating from $p'$ 
(joining $p'$ to some other vertex $q'$, say),
and $e'_2$ is either of the two remaining edges emanating from $q'$,
then there is a unique way to choose $e'_3,\dots$
so that the resulting path shadows $\pi$.
Magically, the signed area enclosed by $\pi'$
is independent of all those choices.

It should also be mentioned that if
we choose basepoints $p$ and $p'$ belonging to $L_1$ rather than $L_0$
then the signed area enclosed by $\pi'$
will be $-I(R)$ rather than $I(R)$.

We will not go into the details of the proof,
except to mention three facts:
the shadow of the counterclockwise boundary of a right-pointing stone
(starting and ending at a boundary point in $L_0$)
encloses signed area 3;
the shadow of the counterclockwise boundary of a left-pointing stone
(starting and ending at a boundary point in $L_0$)
encloses signed area $-3$;
and the shadow of the counterclockwise boundary of a bone
encloses signed area 0.
It would be possible to re-prove the main facts of~\cite{CoLa} and~\cite{Thur}
that pertain to tribone tilings in an elementary fashion 
(without use of combinatorial group theory or combinatorial topology)
along lines laid out in~\cite{MoPa} and~\cite{Pro1}.

\section{Computing the shadow-path area} \label{sec-formula}

In this section we apply the ideas of the previous section
to compute the Conway-Lagarias invariant of the $(a,b)$-benzel.
An extra complication arises from the spurs 
seen in Figures~\ref{fig:corners}~and~\ref{fig:three-cases}.
A spur in a path consists of two consecutive occurrences 
of the same edge (traversed in opposite directions).
Our definition of shadow-paths assumed that $\pi$
was the boundary of a region $R$ and hence did not involve spurs;
we want to allow there to be some $i$'s for which $e_{i} = e_{i+1}$.
Say that such a spur is {\em isolated}
if the shortened path $e_{i-2}, e_{i-1}, e_{i+2}, e_{i+3}$
obtained by removing $e_{i}$ and $e_{i+1}$
satisfies $e_{i-2} \neq e_{i-1} \neq e_{i+2} \neq e_{i+3}$.
Let us suppose that all spurs in $\pi$ are isolated.
We will say that $\pi'$ shadows $\pi$
if $\pi'$ has isolated spurs exactly where $\pi$ does
(and no other spurs at all),
and if for all $i$ such that $e_{i} = e_{i+1}$
the shortened path $e'_{i-2}, e'_{i-1}, e'_{i+2}, e'_{i+3}$
shadows $e_{i-2}, e_{i-1}, e_{i+2}, e_{i+3}$.

\begin{figure}
\begin{center}
\includegraphics[width=3.6in]{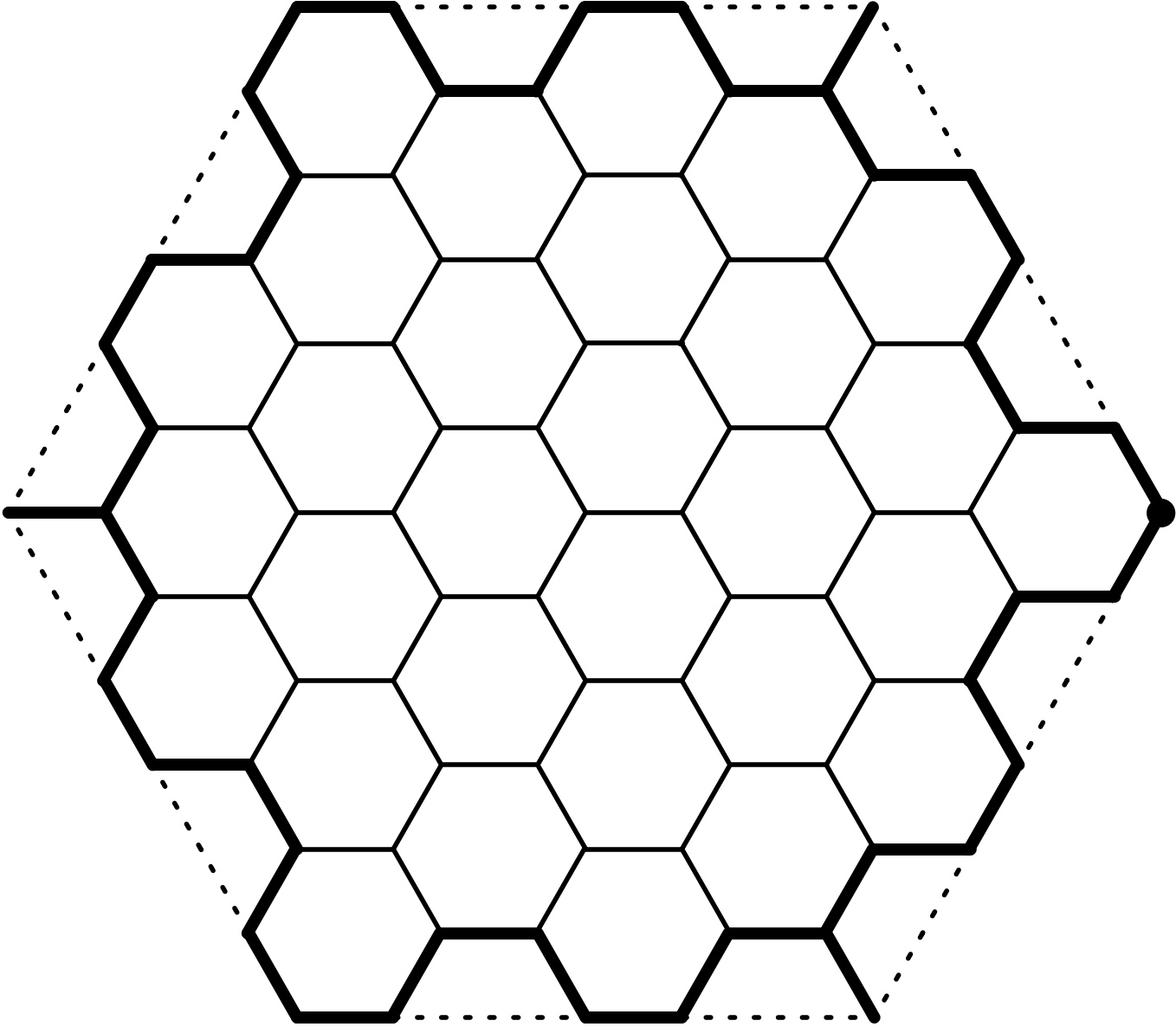}
\end{center}
\caption{The boundary of the (6,6)-benzel.}
\label{fig:six-six}
\end{figure}

\begin{figure}
\begin{center}
\includegraphics[width=3.0in]{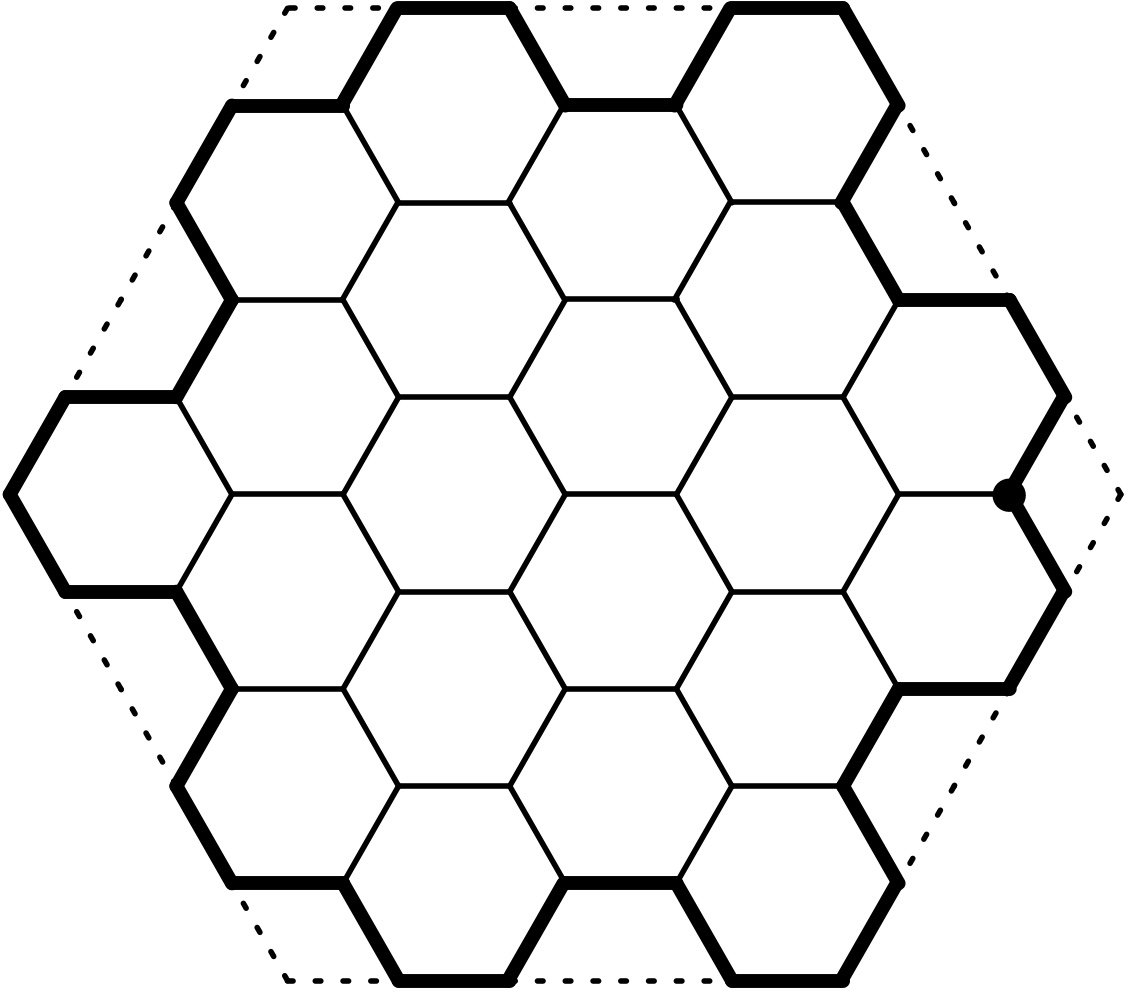}
\end{center}
\caption{The boundary of the (5,5)-benzel.}
\label{fig:five-five}
\end{figure}

\begin{figure}
\begin{center}
\includegraphics[width=2.4in]{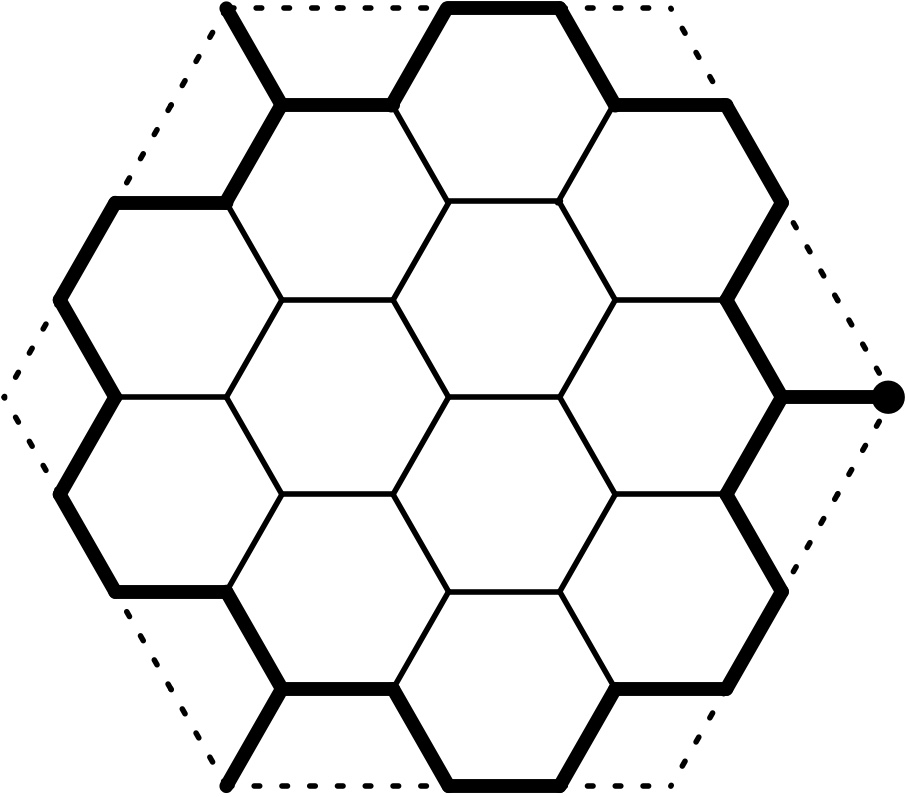}
\end{center}
\caption{The boundary of the (4,4)-benzel.}
\label{fig:four-four}
\end{figure}

We can now describe the boundary words of the benzels and their shadows,
expressed as sequences of unit vectors whose spurs are all isolated.
The reader may find it helpful to consult
Figures \ref{fig:six-six}, \ref{fig:five-five}, and \ref{fig:four-four}
which show the (6,6)-benzel, (5,5)-benzel, and (4,4)-benzel
and are representative of benzels of class 0, 1, and $-1$, respectively.
In the interest of legibility, we omit commas in sequences of unit vectors. 

For class 0 benzels we choose our basepoint $p \in L_0$ to be
the rightmost corner of the bounding hexagon.
Let $s=(2a-b)/3$ and $t=(2b-a)/3$.
The boundary word of the $(a,b)$-benzel is
$$(\bb\ba'\bb\bc')^s (\bc\ba'\bb\ba')^t 
  (\bc\bb'\bc\ba')^s (\ba\bb'\bc\bb')^t 
  (\ba\bc'\ba\bb')^s (\bb\bc'\ba\bc')^t $$
which is shadowed by
$$(\ba\bb'\bc\bb')^s (\bb\ba'\bb\bc')^t  
  (\bb\bc'\ba\bc')^s (\bc\bb'\bc\ba')^t 
  (\bc\ba'\bb\ba')^s (\ba\bc'\ba\bb')^t .$$

For class 1 benzels we choose our basepoint $p \in L_1$ 
so that $p+1$ is the rightmost corner of the bounding hexagon.
Let $s=(2a-b-2)/3$ and $t=(2b-a-2)/3$.
The boundary word of the $(a,b)$-benzel is
$$(\bc'\bb\ba'\bb)^s \bc'\bb (\ba'\bc\ba'\bb)^t \ba'\bc
  (\ba'\bc\bb'\bc)^s \ba'\bc (\bb'\ba\bb'\bc)^t \bb'\ba
  (\bb'\ba\bc'\ba)^s \bb'\ba (\bc'\bb\bc'\ba)^t \bc'\bb$$
which is shadowed by
$$(\ba'\bb\ba'\bc)^s \ba'\bb (\ba'\bb\bc'\bb)^t \ba'\bb 
  (\bc'\ba\bc'\bb)^s \bc'\ba (\bc'\ba\bb'\ba)^t \bc'\ba
  (\bb'\bc\bb'\ba)^s \bb'\bc (\bb'\bc\ba'\bc)^t \bb'\bc .$$

For class $-1$ benzels we choose our basepoint $p \in L_1$ to be
the rightmost corner of the bounding hexagon.
Let $s=(2a-b-1)/3$ and $t=(2b-a-1)/3$.
The boundary word of the $(a,b)$-benzel is
$$(\ba'\bb\bc'\bb)^s \ba' (\bb\ba'\bc\ba')^t \bb
  (\bb'\bc\ba'\bc)^s \bb' (\bc\bb'\ba\bb')^t \bc
  (\bc'\ba\bb'\ba)^s \bc' (\ba\bc'\bb\bc')^t \ba $$
which is shadowed by
$$(\bb'\ba\bb'\bc)^s \bb' (\ba\bc'\ba\bb')^t \ba 
  (\ba'\bc\ba'\bb)^s \ba' (\bc\bb'\bc\ba')^t \bc
  (\bc'\bb\bc'\ba)^s \bc' (\bb\ba'\bb\bc')^t \bb .$$

We do not need to use these shadow words in any explicit way;
all that matters is that, like the boundary words they shadow,
each intermixes fixed strings,
strings raised to the power of $s$,
and strings raised to the power of $t$;
as in section~\ref{sec-area},
we can conclude that the signed area
must be a quadratic function of $s$ and $t$,
which in each of the three cases amounts to
a quadratic function of $a$ and $b$.
(If one wishes to do this explicitly,
one must remember that when $c$ is 1 or $-1$, 
the basepoint $p$ is in $L_1$ rather than $L_0$,
which introduces a sign-change.)

For each $c$ in $\{0,1,-1\}$, one can compute
$I(R)$ for half a dozen small benzels of class $c$.
(The most pleasant way to do this
is to tile $R$ with bones and stones;
then $i(R)$ is the number of right-pointing stones
minus the number of left-pointing stones,
and $I(R)$ is $3i(R)$.)
Equating these values of $I(R)$
with the values taken on by a quadratic with undetermined coefficients,
one obtains a linear system with 6 equations and 6 unknowns,
which one can solve, obtaining the values of the coefficients
and proving the following theorem.


\bigskip

{\bf Theorem 2}: Let $R$ be a benzel of class $c$.  
Then $I(R)$, the (unrescaled) Conway-Lagarias invariant of $R$, is 
$$\left\{ \begin{array}{ll}
(-3a^2+6ab-3b^2+a+b)/2 & \mbox{if $c$ is 0,} \\
(a^2-4ab+b^2+a+b-2)/2 & \mbox{if $c$ is 1,} \\
(-3a^2+6ab-3b^2-a-b+2)/2 & \mbox{if $c$ is $-1$.}
\end{array} \right.$$

\bigskip

\section{When the invariant vanishes} \label{sec-vanish}

We now show that the Conway-Lagarias invariant of the $(a,b)$-benzel vanishes
precisely when $\{a,b\} = \{k(3k+1)/2,k(3k-1)/2\}$ for some positive integer $k$.

The case $a+b \equiv 1$ mod 3 is easily disposed of.
Recall from Theorem 2 that in this case we have
$I(R) = (a^2-4ab+b^2+a+b-2)/2$,
so we are trying to find positive integer 
solutions to $a^2-4ab+b^2+a+b-2 = 0$
satisfying the constraints
$2 \leq a \leq 2b$, $2 \leq b \leq 2a$.
We cannot have $b=2a$ since that would give 
$a+b = 3a \not\equiv 1$,
so $b \leq 2a-1$, and likewise $a \leq 2b-1$; 
but this implies that
$a^2-4ab+b^2+a+b-2 = a(a-2b+1) + b(b-2a+1) - 2 < 0$.

Just as easy to eliminate is the case $a+b \equiv -1$ mod 3. 
Theorem 2 tells us that in this case we have
$I(R) = (-3a^2+6ab-3b^2-a-b+2)/2$,
so we seek positive integers $a,b$
satisfying $-3a^2+6ab-3b^2-a-b+2 = 0$
with $2 \leq a \leq 2b$ and $2 \leq b \leq 2a$.
If we solve for $b$ in terms of $a$
we get a quadratic equation with discriminant 
$(6a-1)^2-12(3a^2+a-2)=25-24a$,
which is negative since $a \geq 2$.

\begin{figure}
\begin{center}
\includegraphics[width=5.0in]{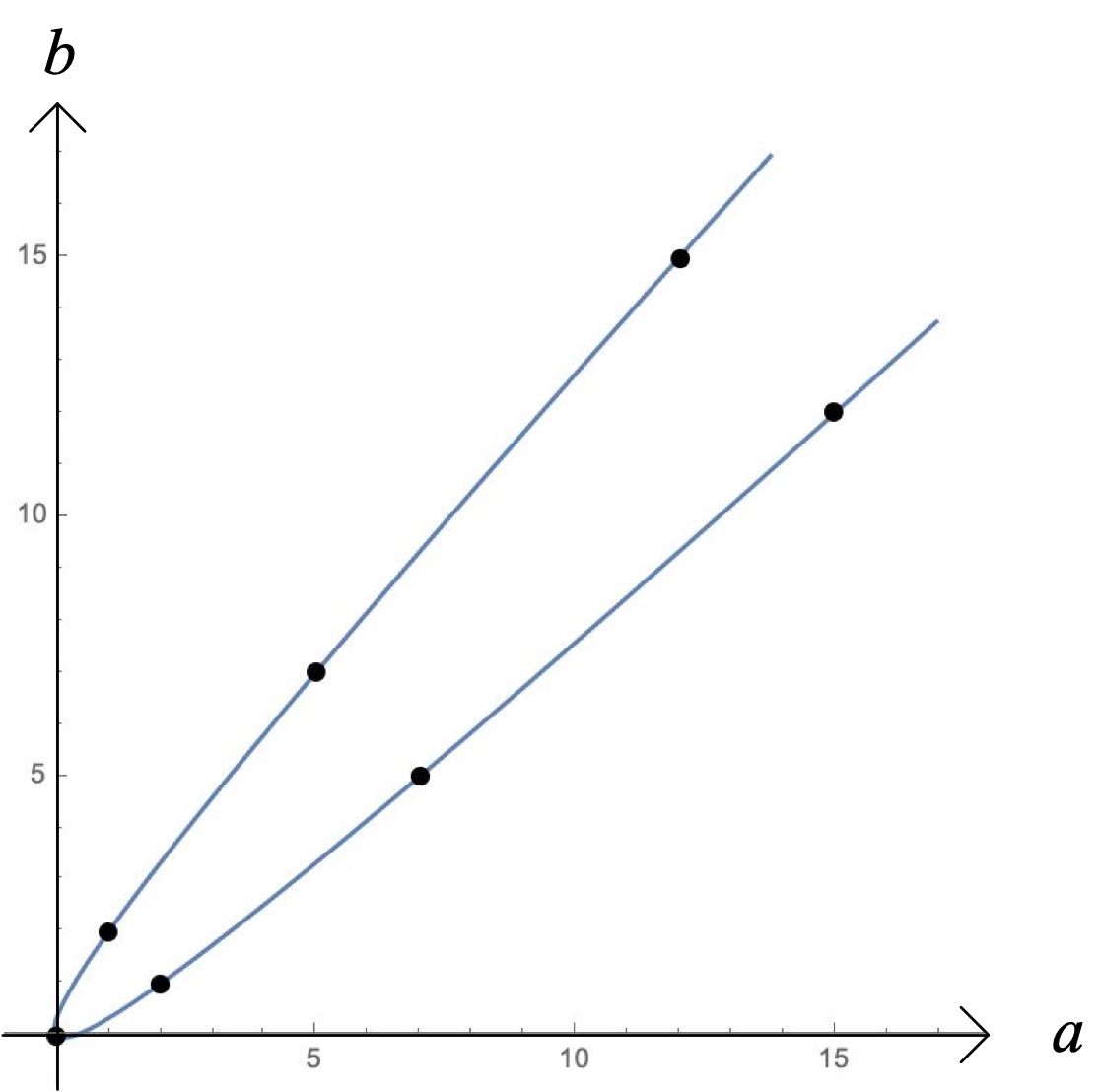}
\end{center}
\caption{Integer solutions to $-3a^2+6ab-3b^2+a+b = 0$.
The solutions (0,0), (1,2), and (2,1) are spurious
for our purposes.}
\label{fig:plot}
\end{figure}

Finally we come to the case $a+b \equiv 0$ mod 3.
Theorem 2 tells us that in this case we have
$I(R) = (-3a^2+6ab-3b^2+a+b)/2$.
Suppose $a$ and $b$ are positive integers 
satisfying $-3a^2+6ab-3b^2+a+b = 0$
with $2 \leq a \leq 2b$ and $2 \leq b \leq 2a$.
If we solve for $b$ in terms of $a$
we get a quadratic equation with discriminant 
$(6a+1)^2-12(3a^2-a)=24a+1$;
this can give us a rational value of $b$
only if $24a+1$ is a perfect square.
Since $24a+1$ is an odd number not divisible by 3,
its integer square roots must be of the form $6k \pm 1$.
Equating $24a+1$ with $(6k+1)^2$ gives $a = (3k^2+k)/2$.
Plugging this into the quadratic gives us an equation
with two solutions, one of which is $b = (9k^2+9k+2)/6$
which is not an integer since
the numerator is not a multiple of 3 and hence not divisible by 6;
the other solution is $b = (3k^2-k)/2$.
Hence we must have $(a,b) = ((3k^2+k)/2,(3k^2-k)/2)$ for some integer $k$.
Likewise, equating $24a+1$ with $(6k-1)^2$ gives 
$(a,b) = ((3k^2-k)/2, (3k^2+k)/2)$.
(Of course both families of solutions are the same,
as can be seen by replacing $k$ by $-k$.)
These solution-pairs are plotted in Figure~\ref{fig:plot};
they lie along a parabola.

\section{Constructing a tiling} \label{sec-construct}

\begin{figure}
\begin{center}
\includegraphics[width=3.0in]{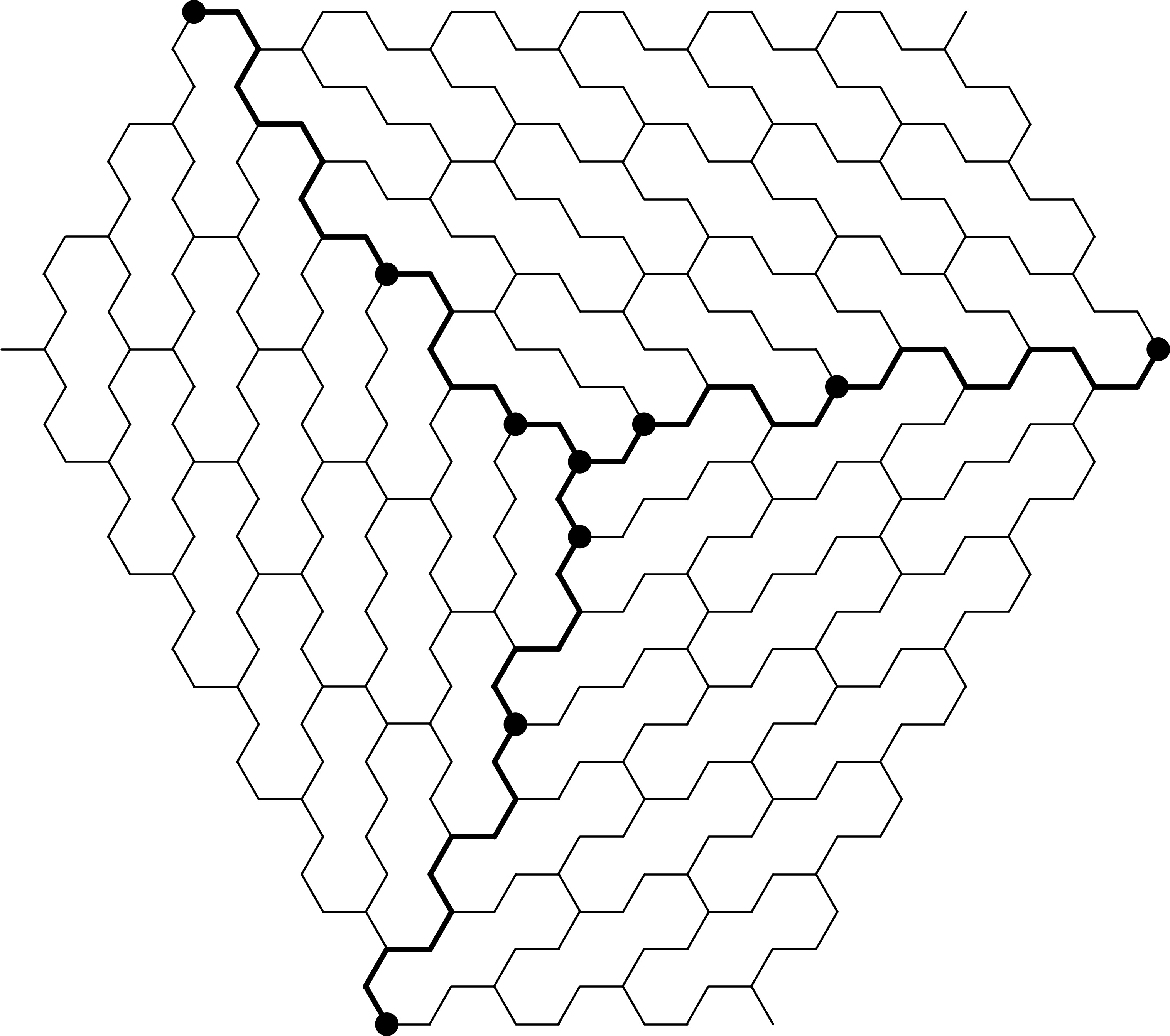}
\end{center}
\caption{Tiling a (12,15)-benzel with bones.
The marked points lie on three parabolas.}
\label{fig:kim}
\end{figure}

In the previous section we gave a condition 
on the ordered pair $(a,b)$ that is necessary 
in order for the $(a,b)$-benzel to admit 
a tiling by tribones; in this section we show 
that the necessary condition is also sufficient.

Figure~\ref{fig:kim} shows a (12,15)-benzel
that has been symmetrically divided into three sectors,
each of which can be tiled by bones in a single orientation.
The path joining the origin to the rightmost point of the benzel
is given by the sequence
$$(\ba\bc'\ba\bb')^0 (\ba\bc') 
  (\ba\bc'\ba\bb')^1 (\ba\bc') 
  (\ba\bc'\ba\bb')^2 (\ba\bc').$$
More generally, consider the path
$$(\ba\bc'\ba\bb')^0 (\ba\bc') 
  (\ba\bc'\ba\bb')^1 (\ba\bc') 
\cdots
  (\ba\bc'\ba\bb')^{k-1} (\ba\bc')$$
We claim that this path and the two paths obtained 
by rotating it 120 and 240 degrees around the origin divide 
the $(a,b)$-benzel with $a=k(3k-1)/2$ and $b=k(3k+1)/2$ into 
three sectors each of which can be tiled by bones in a single orientation.

First, note that the sum of the vectors in the aforementioned path is
\begin{eqnarray*}
&   & \ \ \,(1+3+5+\dots+(2k-1))\ba \\
&   & + (0+1+2+\dots+(k-1))\bb' \\
&   & + (1+2+3+\dots+k)\bc' \\
& = & k^2 \ba - (k(k-1)/2) \bb - (k(k+1)/2) \bc \\
& = & k^2 - (k(k-1)/2) \omega - (k(k+1)/2) \omega^2 \\
& = & - (k^2 + k(k-1)/2) \omega - (k^2 + k(k+1)/2) \omega^2 \\
& = & - ((3k^2-k)/2) \omega - ((3k^2+k)/2) \omega^2 \\
& = & - a \omega - b \omega^2
\end{eqnarray*}
so that the path does indeed lead from
the center of the benzel to the rightmost point
of the bounding hexagon.

\begin{figure}
\begin{center}
\includegraphics[width=4.0in]{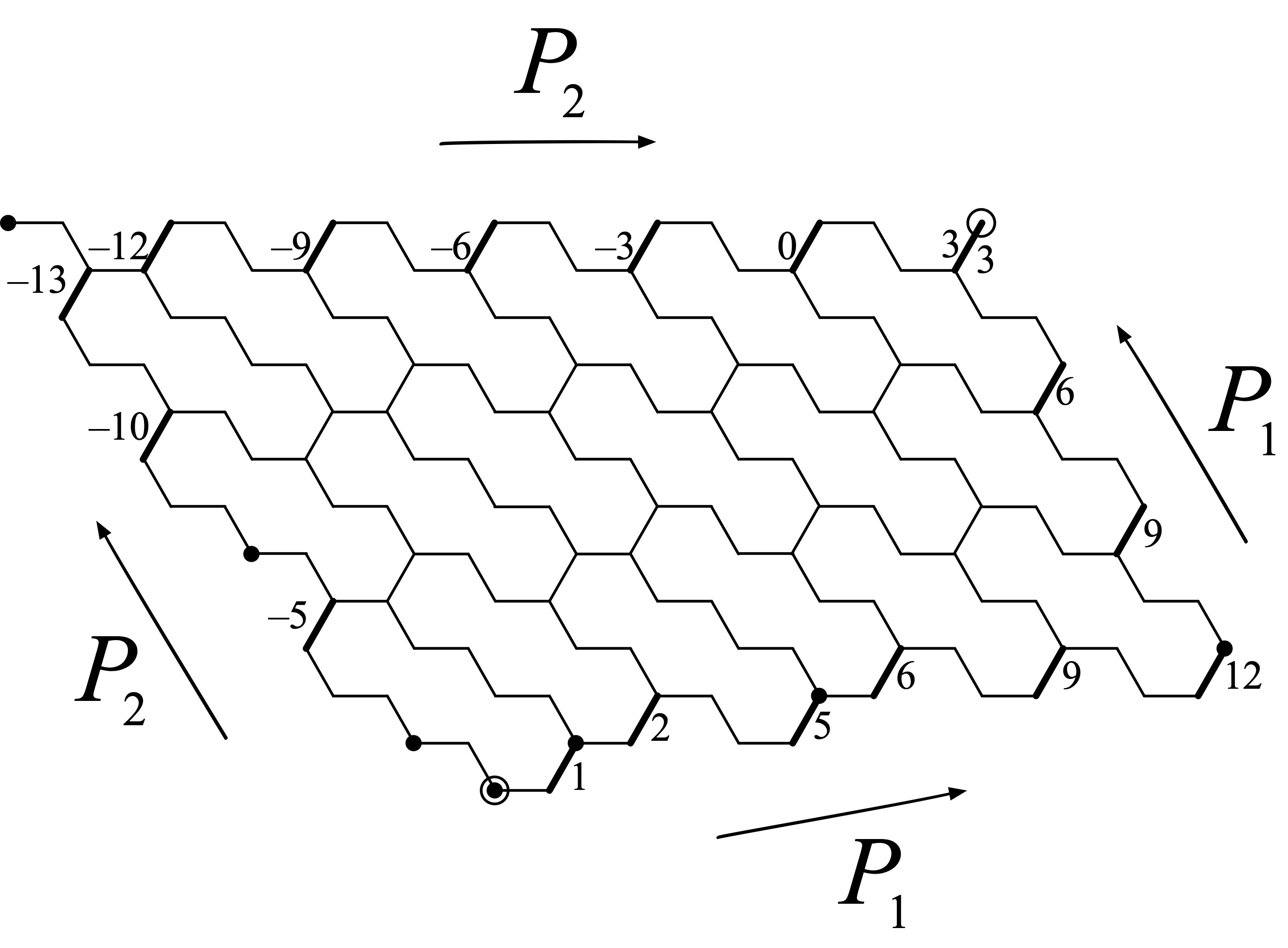}
\end{center}
\caption{One-third of a (12,15)-benzel. The northwest-to-southeast bands, 
progressing from the southwest corner to the northeast corner,
are tiled by 1, 2, 3, 3, 3, 3, 2, 1, and 0 bones, respectively.}
\label{fig:third}
\end{figure}

Put $s = (2a - b)/3 = k(k-1)/2$ and $t = (2b - a)/3 = k(k+1)/2$. 
Consider the path 
$$P_1 =
  (\ba\bc'\ba\bb')^0     (\ba\bc') 
  (\ba\bc'\ba\bb')^1     (\ba\bc') \cdots
  (\ba\bc'\ba\bb')^{k-1} (\ba\bc')
  (\bb\ba'\bb\bc')^{s}$$
joining the origin to the upper right corner of the bounding hexagon
by way of the right corner of the bounding hexagon, and the path 
$$P_2 = 
  (\bb\ba'\bb\bc')^0     (\bb\ba') 
  (\bb\ba'\bb\bc')^1     (\bb\ba') \cdots
  (\bb\ba'\bb\bc')^{k-1} (\bb\ba')
  (\ba\bb'\ba\bc')^{t}$$
joining the origin to the upper right corner of the bounding hexagon
by way of the upper left corner of the bounding hexagon,
as shown in Figure~\ref{fig:third}.
Let $R'$ be the sector of the benzel
contained between the paths $P_1$ and $P_2$.
We need to show that $R'$ can be tiled by translates of
the bone whose axis points in the $\ba\!-\!\bb$ direction.
Assign each point $i\ba+j\bb+k\bc$ ``height'' $i-j$
(readily checked to be well-defined even though $i$, $j$, and $k$ aren't),
so that two points that differ by $\pm\bc$ have the same height;
a bone whose axis points in the $\ba\!-\!\bb$ direction
has two $\pm\bc$-edges on which the height is constant,
and the heights of the two edges differ by 6.
Mark each $\pm\bc$-edge on the boundary of $R'$
with the common height of its two endpoints.
We need to show that the successive heights of the $\pm\bc$-edges
that lie along the path $P_1$ 
($1,2,5,6,9,12,9,6,3$ in the example)
and the successive heights of the $\pm\bc$-edges
that lie along the path $P_2$
($-5,-10,-13,-12,-9,-6,-3,0,3$ in the example)
agree mod 6; that is, we need to check that
subtracting the latter from the former gives a sequence of multiples of 6.
(In the example, subtracting the two sequences and dividing by 6 
gives $1,2,3,3,3,3,2,1,0$;
the terms count the number of $\ba\!-\!\bb$ bones in each diagonal band of $R'$.)
To prove the divisibility claim, it suffices to show that the same claim holds
if we insert a 0 at the beginning of each of the original two sequences
and then replace each resulting sequence by its sequence of first differences;
in the example, the two difference sequences are
$1,1,3,1,3,3,-3,-3,-3$ and $-5,-5,-3,1,3,3,3,3,3$,
and subtracting the second from the first gives $6,6,6,0,0,0,-6,-6,-6$.

The first of the two difference sequences
can be obtained by computing the height-changes
in the subwords between successive occurrences of $\pm\bc$ in $P_1$.
Those $k^2$ subwords yield height-changes
\[
   1, 3^{(0)}, 1, 3^{(1)}, 1, 3^{(2)}, \cdots, 1, 3^{(k-1)}, (-3)^{(s)}
\]
where $h^{(i)}$ denotes $i$ consecutive height changes of size $h$.

The second of the two difference sequences
can be similarly obtained by computing the height-changes
in the subwords between successive occurrences of $\pm\bc$ in $P_2$.
Those $k^2$ subwords yield height-changes
\[
    -5, (-3)^{(0)}, -5, (-3)^{(1)}, -5, (-3)^{(2)}, \cdots; -5, (-3)^{(k-2)}, 1, 3^{(t-1)}.
\]

The difference between the two sequences consists of 
$k(k-1)/2$ 6's, followed by $k$ 0's, followed by $k(k-1)/2$ $-6$'s.
Since all terms are divisible by 6,
we have proved that $R'$ can be tiled by bones
(specifically, bones parallel to $\ba\!-\!\bb$),
which in turn shows that $R$ can be tiled by bones.

Combining the preceding result with the result 
of section~\ref{sec-vanish}, we conclude:

\bigskip

{\bf Theorem 3}: An $(a,b)$-benzel can be tiled by bones
if and only if we have $a = k(3k-1)/2$ and $b = k(3k+1)/2$
or vice versa for some integer $k \geq 2$. 

\bigskip

We have shown that when the Conway-Lagarias invariant of a benzel is zero
the benzel has a tiling that uses only bones.
It appears that when the Conway-Lagarias invariant of a benzel is positive
the benzel has a tiling that uses only bones and right-pointing stones
and that when the Conway-Lagarias invariant of a benzel is negative
the benzel has a tiling that uses only bones and left-pointing stones,
but we do not have a proof.  (See Problems 4 and 7 in~\cite{Pro2}.)
In light of the remarks made at the end of section~\ref{sec-shadow}
we know that these conditions are necessary, 
but we do not know whether they are sufficient.

\section{Further questions} \label{sec-questions}

David desJardins' {\tt TilingCount} program
was useful in leading the second author to conjecture
that bone tilings of the $(a,b)$-benzel
do not exist except when $a$ and $b$
are paired pentagonal numbers.
However, we were not able to get much data
regarding the number of bone tilings that exist
when $a$ and $b$ are of that form,
simply because the area of the $((3k^2-k)/2$,$(3k^2+k)/2)$-benzel
grows like the fourth power of $k$,
and the number of tilings appears to grow 
like an exponential function of the area.
For $k = 2$, the number of tilings is 2; 
for $k = 3$, the number of tilings is $42705 = (3)^2 (5) (13) (73)$;
for $k = 4$, the number of tilings is $7501790059160666750 =
(2) (5)^3 (13) (19) (373) (1559) (208916623)$;
and for $k = 5$, the program runs out of memory.

\begin{figure}
\begin{center}
\includegraphics[width=4.0in]{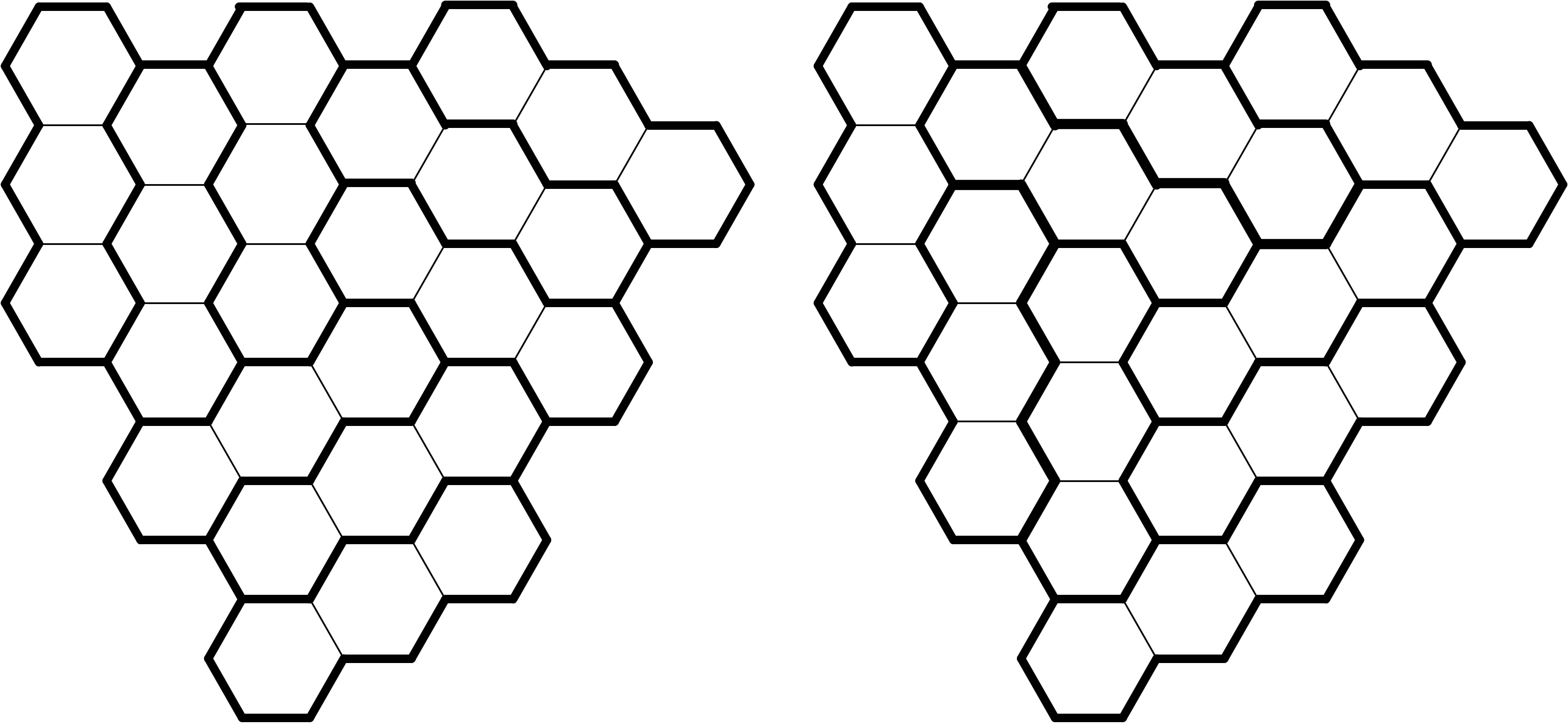}
\end{center}
\caption{The two tilings of the (5,7)-benzel.}
\label{fig:both}
\end{figure}

This article began with a tiling of the (5,7)-benzel.
In fact, there are two such tilings;
Figure~\ref{fig:both} shows both of them.
One can see that the three outermost tiles are actually forced,
in the sense that both tilings contain those tiles.
Similarly, for the (12,15)-benzel,
the tile at the upper right of Figure~\ref{fig:kim} is forced.
It would be good to know more generally
which tiles are forced in a bone tiling of a benzel.

Given that both tilings of the (5,7)-benzel
contain the same number of bones in each of the three orientations,
one might wonder whether this is true for larger benzels. 
In fact this surmise is false;
there exist bone tilings of the (12,15)-benzel
with different numbers of bones in the three orientations.

Finally, one finds that some tiles near the boundary of the (12,15)-benzel,
although not actually forced, are highly likely to occur.
For instance, consider the bone in Figure~\ref{fig:kim}
immediately to the left of the aforementioned forced tile;
it occurs in 42587 of the 42705 tilings,
which means that one sees it more than 99.7 percent of the time
if one picks a tiling uniformly at random.
This suggests that there is a (possibly large) vicinity of the boundary
that is statistically ``frozen'' in the sense that 
a uniformly random tiling exhibits very little variety there.
Is this frozen region (assuming it exists) a thin layer
near the boundary of the benzel,
or does it penetrate a macroscopic distance into the interior?

\subsection*{Acknowledgements}
The authors thank David desJardins, without whose code
this main result would never have come to light;
Timothy Chow, who suggested the term ``phone'' 
to refer to the third trihex tile;
Jonathan Boretsky, 
Son Nguyen, 
and Ashley Tharp 
who along with the first author explored bone tilings of
$(k(3k-1)/2,k(3k+1)/2)$-benzels;
Colin Defant and Rupert Li,
who gave helpful comments on earlier versions of this writeup;
Juri Kirillov, who made some useful observations;
and the two anonymous referees for their careful reading of the manuscript
and their helpful suggestions.
Thanks also to the organizers of the 
Open Problems in Algebraic Combinatorics 2022 conference
for inviting the second author to speak and for funding his participation;
the second author did not know how to prove the main result of section 7 
until students attending the conference, including the first author,
figured it out.

\end{document}